\documentclass[12pt]{article}
\usepackage[T1]{fontenc}
\usepackage{amssymb,amsmath}
\usepackage{epsfig,color}
\usepackage{url,hyperref}
\usepackage[all,dvips,import]{xy}

\title{Compact generation for Lie groupoids}
\author{\small Nicolas Raimbaud}
\date{\today}

\definecolor{myred}{RGB}{247,79,56}
\definecolor{myblu}{RGB}{96,185,236}
\definecolor{lgray}{gray}{.7}
\newxyColor{lagray}{.7}{gray}{}

\def\gc{\color{lgray}}
\def\ds{\displaystyle}
\def\df{\bf}

\def\bc{\begin{center}}
\def\ec{\end{center}}
\def\nin{\noindent}
\def\m#1{\mbox{#1}}

\DeclareMathOperator*{\cups}{\cup}

\newcounter{main}[section]
\renewcommand{\themain}{\thesection.\arabic{main}}
\def\main{\refstepcounter{main}\nin{\bf\themain} }
\def\notmain{\refstepcounter{main}\nin{\themain} }

\newenvironment{basestyle}[2]
{\vspace{0.6\medskipamount} #1\\[.5ex]
\begin{tabular}{cl}#2
\begin{minipage}{.93\textwidth}\em}
{\end{minipage}
\end{tabular}\vspace{1.4\medskipamount} \\}
\newenvironment{thmstyle}[1]
{\begin{basestyle}{#1}{\hspace{1pt}{\gc\vrule width 1pt}&\hspace{-.8em}}}
{\end{basestyle}}
\newenvironment{defstyle}[1]
{\begin{basestyle}{#1}{{\gc\vrule width 1pt}\hspace{1pt}{\gc\vrule width 1pt}&\hspace{-.7em}}}
{\end{basestyle}}
\newenvironment{lemma}[1]
{\begin{thmstyle}{\\ \sc Lemma \notmain #1}}
{\end{thmstyle}}
\newenvironment{prop}[1]
{\begin{thmstyle}{\\ \sc Proposition \notmain #1}}
{\end{thmstyle}}
\newenvironment{mainprop}[1]
{\begin{thmstyle}{\main\sc Proposition #1}}
{\end{thmstyle}}
\newenvironment{thm}[1]
{\begin{thmstyle}{\\ \sc Theorem \notmain #1}}
{\end{thmstyle}}
\newenvironment{mainthm}[1]
{\begin{thmstyle}{\main\sc Theorem #1}}
{\end{thmstyle}}
\newenvironment{cor}[1]
{\begin{thmstyle}{\\ \sc Corollary \notmain #1}}
{\end{thmstyle}}

\newenvironment{defenv}[1]
{\begin{defstyle}{\\ \sc Definition \notmain #1}}
{\end{defstyle}}
\newenvironment{maindef}[1]
{\begin{defstyle}{\main\sc Definition #1}}
{\end{defstyle}}

\def\proof{{\em Proof.} }
\def\cqfd{$\Box$}
\def\eq{\;=\;}
\newcommand{\f}[3]{\ensuremath{#1:#2 \to #3}}
\def\R{\mathbb{R}}
\def\C{\mathbb{C}}
\def\T{{\rm T}}
\def\id{{\rm id}}

\def\gdim{{\rm gdim}\,}
\def\tloc{{\mathcal T}}
\def\Hol{Hol}
\def\Mon{Mon}
\def\Orb{Orb}
\def\fol{{\cal F}}
\def\pgH{{\cal H}}
\def\gK{{\cal K}}

\def\gU{{\cal U}}
\def\bgU{\gU_0}
\def\gV{{\cal V}}
\def\bgV{\gV_0}
\def\gC{{\cal C}}
\renewcommand{\phi}{\varphi}

\renewcommand{\theta}{\vartheta}
\def\Db{D_{m^*}}
\def\conj{\overline{\phantom{:}.\phantom{:}}}

\begin{document}
\maketitle

\begin{abstract}
 Thirty years after the birth of foliations in the 1950's, Andr\'e Haefliger has introduced a special property satisfied by holonomy pseudogroups of foliations on compact manifolds, called compact generation. Up to now, this is the only general property known about holonomy on compact manifolds.\par
 In this article, we give a Morita-invariant generalization of Haefliger's compact generation, from pseudogroups to object-separated Lie groupoids.
\end{abstract}

\nin {\em Mathematics Subject Classification:} 57R30, 58H05.\\
\nin {\em Keywords:} Foliation, Lie groupoid, Morita equivalence, Compact generation.

\section*{Introduction}

Recall that a foliation of an $n$-manifold $M$ is a partition of $M$ into $p$-submanifolds, which locally looks like parallel copies of $\R^p$ in $\R^n$. The connected components of the $p$-structure are called the leaves of the foliation.\\
A classical tool to study the dynamics of the leaves is the set of holonomy elements of the foliation, which are diffeomorphisms between transversals (submanifolds everywhere transverse to the leaves, and of complementary dimension). Given a complete transversal $T$ (a transversal that meets every leaf at least once), the holonomy elements between open sets of $T$ generate a pseudogroup, the holonomy pseudogroup of the foliation relative to $T$. This pseudogroup essentially depends on the foliation, in the sense that two different complete transversals will provide two Haefliger-equivalent [\ref{pg2gd}] pseudogroups. When the base manifold $M$ is compact, Haefliger introduced an invariant `finite generation' property for the holonomy pseudogroups: compact generation.\\
Given a pseudogroup $(\pgH, T)$, its set $H$ of germs is naturally equipped with a composition operation (the compositoin of germs). This composition makes $H$ a small category with set of objects $T$, in which every arrow is invertible. Such a structure is called a groupoid structure. The set of germs $H$ can also be naturally given a topology (the sheaf topology), and natural charts with smooth changes of coordinates inherited from $T$. This differential structure is compatible with the algebraic structure of $H$, and eventually $H$ is a very special kind of groupoid, an (effective) \'etale Lie groupoid.\\
In \cite{Haef00}, Andr\'e Haefliger translated the compact generation property for pseudogroups in terms of \'etale groupoids, using the close relation between pseudogroups of diffeomorphisms and their groupoids of germs [\ref{pg2gd},\ref{gd2pg}]. In this paper, we extend compact generation to (almost) all Lie groupoids, in a consistent way: invariance of compact generation under Haefliger equivalence is extended to invariance under Morita equivalence. I would like here to thank Ga\"el Meigniez for leading me to the right definition through our many dicussions. I still can't realize the amount of `groupoidic' material we studied together before the final formulation arose.

\section{Preliminaries}

\main Recall that a {\df groupoid} $G$ is a small category in which every arrow is invertible. Its set of objects, classically denoted $G_0$, is called the {\df base} of the groupoid, and can be embedded in $G$ through the {\df unit map} (sending each object to its associated identity morphism). We shall always use this identification, so that $G_0\subset G$ for any groupoid. When needed, the set of morphisms is distinguished from the groupoid $G$ by denoting it $G_1$. The canonical maps $G_1\to G_0$ sending each morphism to its starting or ending object are called the {\df source and target maps} of $G$ and denoted $s$ and $t$, respectively.\\
A groupoid $G$ may be seen as a set of {\df points} $G_0$, together with {\df arrows} between these points, where the composition of morphisms is some kind of concatenation of arrows. We shall use a leftward-pointing convention for groupoid arrows: an arrow (morphism) $g$ from $x$ to $y$ shall be written \f{g}{x}{y}, and if \f{g}{x}{y} and \f{h}{y}{z} are two arrows of $G$, the composition of $g$ and $h$ will be \f{hg}{x}{z}.\\
The basic algebraic behavior of groupoids is very similar to groups' \cite{MackenzieI}, therefore both the vocabularies of categories and groups is usually applied to groupoids: the composition is called ``{\df product}'', identity morphisms are called ``units'', and so on. In this article, we shall be especially interested in the obvious notions of subgroupoid (a subset closed under composition and inversion), and of {\df subgroupoid generated by a set} (the intersection of all subgroupoids containing that set). We will also call ``full'' a subgroupoid if it is full as a subcategory.

\nin {\bf Notations:} If $X$ and $Y$ are two subsets of $G_0$, the set of all arrows issuing from any point of $X$ and ending at any point of $Y$ will be denoted $G_X^Y$ (if $X$ or $Y$ is the whole base $G_0$, we won't specify it, writing $G^Y$ or $G_X$).

\begin{maindef}{[Lie groupoid]}
 A groupoid $G$ is called a {\df Lie groupoid} if $G$ and $G_0$ are smooth finite-dimensional {\em not necessarily Hausdorff} manifolds, if the source and target maps of $G$ are submersions, and if its product, inverse and unit maps are smooth.
\end{maindef}
By assuming $s$ and $t$ to be submersions, the fibered product $G_s\times_t G=G\times_{s=t}G$ inherits a smooth submanifold structure, and that's why smoothness of the product makes sense. Also note that $G_0\subset G$ is a submanifold, as image of a section (unit map) of a submersion (source or target map).\\
It is important here to allow $G$ to be non-Hausdorff ; it is indeed often the case in practice with groupoids arising in foliation theory.
We also allow non-Hausdorff bases, but for convenience only, so that we don't have to check the Hausdorff condition in our manipulations (however, we will remove this liberty on bases in the last section). Anyway, from now on, {\em no manifold is considered Hausdorff unless explicitly stated}\footnote{It may be useful to recall that in any case, points in a topological manifold are always closed (the $T_1$ axiom is satisfied).}. This lack of Hausdorff condition does not change many things, as long as you remember never to use closures of sets (which may be pretty wild).\label{relk} Especially, as compact sets are not always closed (that is, are only quasi-compact in the Bourbaki sense), {\df relative compactness} shall be understood as `being included in a compact set', which is a (strictly) weaker condition than `having a compact closure'.

\main In group theory, any element of a group naturally generates two diffeomorphisms of the group, the left and right translations by this element. There's no immediate equivalent in groupoid theory: given a point $x$ in the base, we may `translate' it with any arrow $g$ issuing from $x$, but there's no way to translate other elements of the base with $g$, for it has only one source: $x$ ! Thus to define a translation on $G$ we at least have to chose for each point $x$ of the base an element starting at $x$. To carry out such choices, we follow Kirill Mackenzie \cite{MackenzieI} and introduce the notion of bisection.
\begin{defenv}{[Bisection]}
 Let $G$ be a Lie groupoid. Any (smooth) section $\beta$ of the source map $s$ of $G$ such that the composition $t\beta$ is a diffeomorphism of $G_0$ is called a global bisection of $G$.
\end{defenv}
Note that this notion is symmetric in $s$ and $t$: when we identify any $s$-section to its image in $G$, bisections are those submanifolds of $G$ for which the restrictions of both $s$ and $t$ are diffeomorphisms onto $G_0$. Bisections are very efficient in groupoid theory, due to the following fact:
\begin{prop}{}\label{lbeta}
 For any element $g$ of a Lie groupoid $G$, $sg\mapsto g$ may be extended to a local {\df bisection} $\beta:U\to G$ (where $U$ is a neighborhood of $sg$).
\end{prop}
A local bisection may be seen, of course, as a smooth submanifold of $G$ for which both $s$ and $t$ are embeddings into $G_0$. This point of view almost yields the proof of the proposition: any small enough smooth submanifold containing $g$, simultaneously transverse to the $s$- and $t$-fibers through $g$ fits. As we have plenty of them, we shall for convenience refer to local bisections as bisections.\\
Given a bisection $\beta:U\to G$, one can easily define a (local left-) {\df translation} $L_{\beta}:G^U\to G^{t\beta(U)}$ by letting any element of the image of $\beta$ act by left-translation on the $t$-fiber of its source: $L_{\beta}(h)=\beta(th)\cdot h$. Such a translation is a diffeomorphism of $G$ (compose $\beta$ with the inversion and build an inverse mapping). To illustrate these notions, let us prove the following:
\begin{prop}{}\label{openprod}
 The product of two open sets in a Lie groupoid is an open set.
\end{prop}
\proof We call $\gU$ and $\gV$ these open sets. Let $u\in\gU$ and $v\in\gV$ with $sv=tu$, we wish to prove that $\gV\cdot\gU$ is a neighborhood of $vu$. Take a bisection $\beta:W\to G$ extending $v$, shrinking it if necessary we may assume $\beta(W)\subset\gV$. Now $L_{\beta}(\gU\cap G^W)$ is an open set, containing $vu$ and included in $\gV\cdot\gU$. \cqfd

\main\label{pbgd} Two different categories of Lie groupoids are usually being used. The straightforward algebraic one:
\begin{defenv}{}
 A {\df Lie groupoid morphism} is a smooth functor between Lie groupoids.
\end{defenv}
\ldots and a more mysterious, intricate, but more significant one, the Hilsum-Skandalis category \cite{MrcunKT,MrcunPhD}. We won't go into too much detail in the Hilsum-Skandalis category, however we need the equivalence relation associated to this category: Morita equivalence. This equivalence is based on the notion of pullbacks in the ``algebraic'' Lie groupoid category, which we introduce now.\\
Consider $H$ a Lie groupoid, and $f:M\to H_0$ any smooth map. As the source of $H$ is a submersion, we can construct the smooth fibered product $H_s\times_f M$.
$$\xymatrix @dr @M=+1mm @!=10mm {
H_s\times_f M
 \ar @[|<1pt>] [r] ^{\ds \pi_M} _{~}="back"
 \ar [d] ^{\pi_H}
& M
 \ar [d] ^f\\
H
 \ar [d] ^t
 \ar [r] _s ^{~}="pull"
&H_0 \\
H_0
 \ar @/_2pc/ @[|<1pt>] _{\ds t\pi_H} "1,1";"3,1"
 \ar @{-->} |{\txt{\scriptsize pullback}} "pull";"back"
}$$
If we consider $H$ as a set of arrows from $H_0$ to $H_0$, we've actually just pulled their sources with $f$ from $H_0$ to $M$, so we may see $H_s\times_f M$ as a set of arrows from $M$ to $H_0$. If $t\pi_H$ is again a submersion, we may build the product $M_f\times_{t\pi_H}(H_s\times_f M)$, that is, we can also pull the targets of the `arrows' of $H_s\times_f M$ with $f$, and enventually get `arrows' from $M$ to $M$.\\
$$\xymatrix @dr @!=13mm {
M_f\times_{t\pi_H}(H_s\times_f M)
 \ar @[|<1pt>] [r] ^{\pi_{HM}}
 \ar @[lagray] [d]
 \ar @/_2pc/ @[|<1pt>] [dd] ^(.3){~}="backa" ^(.72){~}="backb" ^{~}="backc"
& H_s\times_f M
 \ar @[|<1pt>] [r] ^{\pi_M} _{~}="back"
 \ar @[lagray] [d] ^{\gc \pi_H}
& M
 \ar @[lagray] [d] ^{\gc f}  \\
{\gc M_f\times_t H}
 \ar @[lagray] [r]
 \ar @[lagray] [d]
&{\gc H}
 \ar @[lagray] [r] _{\gc s} ^{~}="pull"
 \ar @[lagray] [d] ^{\gc t}
&{\gc H_0}  \\
M
 \ar [r] _f
&{H_0}
\ar @/_2pc/ "1,2";"3,2" _(.52){t\pi_H} _(.3){~}="pulla" _(.7){~}="pullb" _{~}="pullc"
\ar @[lagray] @{-->} "pull";"back"
\ar @{-->} "pulla";"backa"
\ar @{-->} "pullb";"backb"
\ar @{} |<(.4){\txt{pullback}} "pullc";"backc"
}$$
It is an elementary exercise to check that under this assumption, the triple product $M_f\times_t H_s\times_f M$ is canonically isomorphic to its two decompositions in nested double fiber products. This triple product has a natural Lie groupoid structure: take the obvious source and target maps, and define a compatible product by $(p,h,n)(n,g,m)=(p,hg,m)$. All structure maps are clearly smooth, and it is easy to check that the source and target maps are submersions. In this situation, $M_f\times_t H_s\times_f M$ is called the {\df pullback groupoid} of $H$ by $f$, and denoted $f^*H$. Note that in particular, this construction can be achieved entirely (i.e. $t\pi_H$ is a submersion) when $f$ is a submersion.

\main\label{ee} Let $\phi:G\to H$ be a Lie groupoid morphism. Then the induced map $\phi_0:G_0\mapsto H_0$ (the base of $\phi$) is a particular case of a map $f:M\to G_0$ considered in the previous paragraph. If $t\pi_H:H_s\times_{\phi_0}G_0)\to H_0$ is a surjective submersion, we may consider the pullback $\phi_0^*H$. It is endowed with a natural Lie groupoid morphism \m{$(t,\phi,s):G\to\phi_0^*H$}. When this morphism is an isomorphism, i.e. when $G$ naturally identifies to $\phi_0^*H$, $\phi$ is called an {\df essential equivalence} \cite{MoerdijkFolGd,MrcunKT}. Note that in particular $(t,\phi,s)$ is a bijection, which means that any arrow \f{h}{\phi_0(x)}{\phi_0(y)} of $H$ between points in the image of $\phi_0$ admits a unique lifting \f{g}{x}{y} in $G_x^y$ ({\df unique lifting property}).
In the particular case when $\phi_0$ is already a surjective submersion, an essential equivalence is called a {\df Morita morphism}.\\
Be aware that existence of an essential equivalence between two groupoids is {\em not} an equivalence relation, for it is not symmetric. Actually, Morita equivalence is the associated symmetrized relation:
\begin{defenv}{[Morita equivalence]}\label{moree}
 Two groupoids $G$ and $H$ are said to be {\df Morita-equivalent} if there exists a third groupoid $K$ and two essential equivalences $\phi,\psi:K\to G,H$.
\end{defenv}
It can be seen that this relation is an equivalence relation. Morita equivalence is a very flexible notion and admits many other equivalent definitions. In particular we shall be interested in the following one (see also \ref{mordim}):
\begin{thm}{}\label{mormor}
 Two Lie groupoids $G$ and $H$ are Morita-equivalent if, and only if, there exists a third groupoid $K$ and two Morita morphisms $\phi,\psi:K\to G,H$.
\end{thm}
This means we can freely assume that the functors used in the definition are submersions on the bases (which is actually equivalent to being a submersion for an essential equivalence).

\section{Groupoidizing pseudogroups}

\main\label{pg2gd} Recall that a {\df pseudogroup} of diffeomorphisms on a manifold $T$ is a set of diffeomophisms between open sets of $T$, which is closed under restriction, inversion and gluing. To any pseudogroup $(\pgH,T)$, one can associate the set $H$ of all germs of elements of $\pgH$, which is a groupoid for the composition of germs $[g]_y\cdot[f]_x=[g\circ f]_x$. This groupoid has a natural (sheaf) topology: given an element $h\in\pgH$, the collection of all germs of $h$ at all points of its domain represents a base open set for this topology. Each such base open set may be identified to an open set of $T$, therefore the changes of coordinates in $H$ identify to changes of coordinates in $T$, and thus are smooth. Moreover, it is obvious in these particular coordinates that the source and target maps of $H$ are \'etale (i.e. local diffeomorphisms), and that its other structure maps are smooth. We shall call $H$ the {\df germ groupoid} of $\pgH$, and denote it $[\pgH]$. It is an {\df \'etale groupoid}, a Lie groupoid with \'etale source and target.\\
Recall the definition of {\df Haefliger equivalence} between pseudogroups (see \cite{Haef00} for another definition): two pseudogroups $(\pgH,T)$ and $(\pgH',T')$ are said to be Haefliger equivalent if there exists a set $\Phi$ of diffeomorphisms from open sets covering $T$ to open sets covering $T'$ such that:
$$\Phi\pgH\Phi^{-1}\ \subset\ \pgH' \qquad\text{and}\qquad \Phi^{-1}\pgH'\Phi\ \subset\ \pgH$$
\vspace*{-1.2em}
\begin{prop}{}\label{haef2mor}
 If two pseudogroups $(\pgH,T)$ and $(\pgH',T')$ are Haefliger-equivalent, their germ groupoids $H$ and $H'$ are Morita-equivalent [\ref{moree}].
\end{prop}
\proof Let $\Phi$ be a set of diffeomorphisms giving an Haefliger equivalence from $\pgH$ to $\pgH'$, and let's write $Z_0=H'[\Phi]H$ the collection of all germs coming from compositions of maps of $\pgH$, $\Phi$ and $\pgH'$ (whenever defined).\\
\centerline{
\xy
 \xyimport(10,10){
  \resizebox{!}{60mm}{ \includegraphics{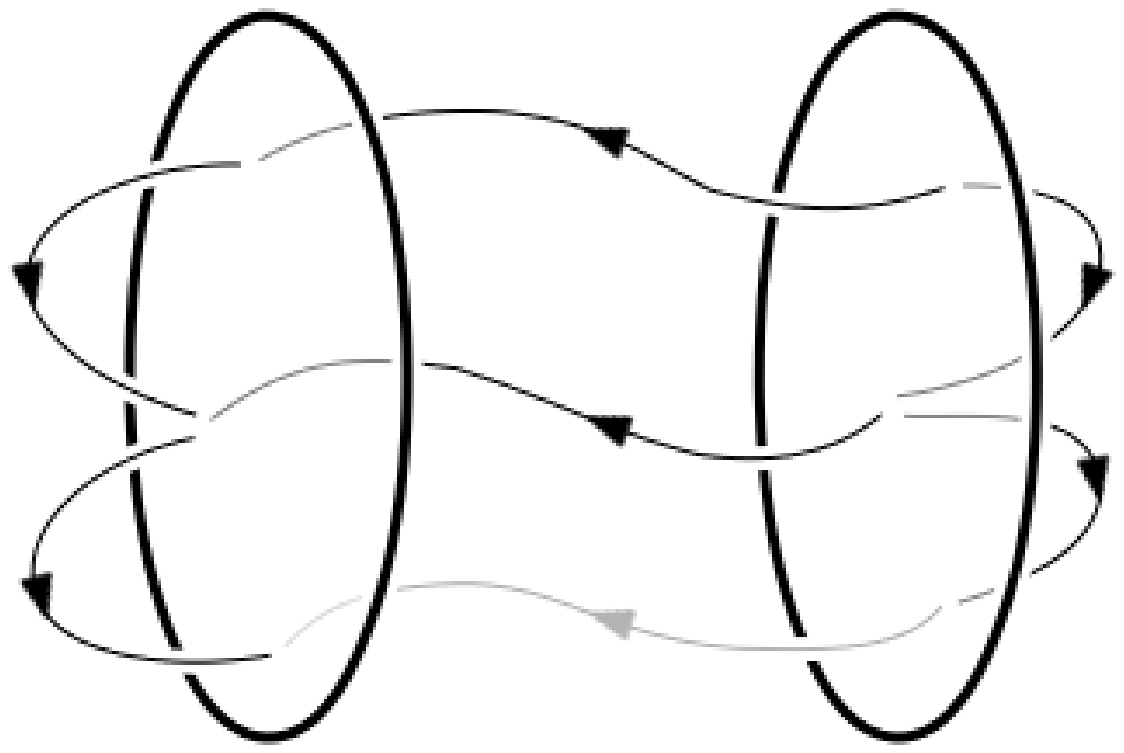} }
  },
 (2.9,.7)*\txt{$T'$},
 (7,.7)*\txt{$T$},
 (5,8)*\txt{$z_1$},
 (5.1,5.4)*\txt{$z_2$},
 (5.1,1.9)*\txt{\gc $h'_2 z_2 h_2^{-1}$},
 (.8,6)*\txt{$h'_1$},
 (.8,3.1)*\txt{$h'_2$},
 (8.8,6.4)*\txt{$h_1$},
 (8.8,4.1)*\txt{$h_2$},
 (5,9.4)*\txt{$Z_0$},
 (0,4.7)*\txt{$H'$},
 (9.5,5.1)*\txt{$H$},
 (6.2,9.4)*\txt{$\sigma$},
 (3.8,9.4)*\txt{$\tau$},
 \ar (5.5,9.3) ; (6.6,8.8),
 \ar (4.5,9.3) ; (3.3,8.8)
\endxy
}
Write $\sigma$ and $\tau$ the maps sending each element of $Z_0$ repectively to its source (in $H$) and its target (in $H'$), and give the manifold $Z=H'_s\times_{\tau} Z_0\, {}_{\sigma}\times_s H$ a source map $s=\pi_{Z_0}$, a target map $t(h',z,h)=h'zh^{-1}$, and a product
$$(h'_2,z_2,h_2)(h'_1,z_1,h_1)=(h'_2 h'_1,z_1,h_2 h_1)$$
It is easy to check that $Z$ is an \'etale groupoid (all maps here are \'etale). It is no more difficult to see that the two obvious smooth functors ${\mathfrak S}=\pi_H:Z\to H$ and ${\mathfrak T}=\pi_{H'}:Z\to H'$ are actually Morita morphisms: apply the pullback construction [\ref{pbgd}] to the maps $\sigma={\mathfrak S}_0$ and $\tau={\mathfrak T}_0$ (\'etale and thus submersive), and use the global composition of germs over $T\cup T'$ to build smooth inverses for the canonical functors (for example with $\mathfrak S$):\vspace*{-0.8em}
$$\xymatrix @W=30mm @C=19mm @R=0pt {
H'_s\times_{\tau} Z_0\, {}_{\sigma}\times_s H
 \ar [r] ^{\cong}
&Z_0\, {}_{\sigma}\times_t H_s \times_{\sigma} Z_0 \\
(h',z_1,h)
 \ar@{|->} [r] ^{(t,{\mathfrak S},s)}
& (h'z_1 h^{-1},h,z_1) \\
(z_2 h z_1^{-1}, z_1, h)
& (z_2,h,z_1)
 \ar@{|->} [l]
}$$
\vspace*{-0.2em} This implies that $H$ and $H'$ are Morita-equivalent. \cqfd

\pagebreak
\begin{maindef}{}\label{gd2pg}
 We shall call {\df $\mathbf 0$-translation} of a Lie groupoid $H$ any diffeomorphism between open sets of $H_0$ which may be locally written $t\beta$, for some local bisection $\beta$.
\end{maindef}
It is easy to check that the $0$-translations of a groupoid $H$ form a pseudogroup of diffeomorphisms on $H_0$. We shall denote this pseudogroup $\tloc H$. In the particular case when $H$ is the germ groupoid of a pseudogroup $\pgH$, the topology of $H$ forces the bisections to be locally written $x\mapsto [f]_x$ for some $f\in\pgH$. Therefore a $0$-translation $\tau$ of $H$ may be in turn locally written $\tau(x)=t([f]_x)=f(x)$, so that $\tau\in\pgH$ locally, and then globally by gluing ($\pgH$ pseudogroup). Thus $\tloc H\subset\pgH$, and the other inclusion is immediate with the identity $f=x\mapsto t([f]_x)$. Finally we see that $\pgH=\tloc H$, and also $H=[\tloc H]$, i.e. $H$ is effective\footnote{See \cite{MrcunPhD} for more information about the `effect-functor' $[\tloc\cdot]$}. (Note that the two operations $[\,\cdot\,]$ and $\tloc$ might be used to identify pseudogroups and effective \'etale groupoids.)\\
As Haefliger equivalence of two pseudogroups implies Morita equivalence of their germ groupoids, it is natural to ask whether it is true in the other direction. Of course, this question makes sense only if we consider groupoids with bases of the same dimension.
\begin{prop}{}\label{mor2haef}
 If $H$ and $H'$ are two Morita-equivalent groupoids with $\dim H_0=\dim H'_0$, then their pseudogroups of $0$-translations are Haefliger-equivalent.
\end{prop}
We still need some technical results about Morita equivalence to prove this proposition, therefore we will postpone the proof to paragraph \ref{m2hproof}. We may sum up the results of the two last paragraphs in the following theorem:
\begin{thm}{}\label{thmpggd}
 Let $(\pgH,T)$ and $(\pgH',T')$ be two pseudogroups. We may identify these pseudogroups to their germ groupoids, for $\tloc[\pgH]=\pgH$ (same for $\pgH'$).\\
 The pseudogroups $\pgH$ and $\pgH'$ are Haefliger-equivalent if and only if the groupoids $[\pgH]$ and $[\pgH']$ are Morita-equivalent.
\end{thm}

\begin{mainprop}{}\label{mordim}
 If $H$ and $H'$ are two Morita-equivalent groupoids with $\dim H_0=\dim H'_0=d$, there exists a third groupoid $Z$ with $\dim Z_0=d$ and two Morita morphisms $\phi,\psi:Z\to H,H'$.
\end{mainprop}
\proof Take $K$ a Lie groupoid yielding a Morita equivalence between $H$ and $H'$, with Morita morphisms $\eta,\theta:K\to H,H'$ (theorem \ref{mormor}). Take any $m\in K_0$, and write $x=\eta_0 m$, $x'=\theta_0 m$. The subspaces $\T_m(\eta_0^{-1}x)$ and $\T_m(\theta_0^{-1}x')$ have the same codimension $d$ in $\T_m K_0$, so they admit a common supplementary $F_m$. Take some coordinates system around $m$, and consider a small $d$-disc $D_m$ containing $m$, and contained in the subspace $F_m$ in those coordinates.
\centerline{
\xy
\xyimport(10,10){
 \resizebox{!}{49mm}{ \includegraphics{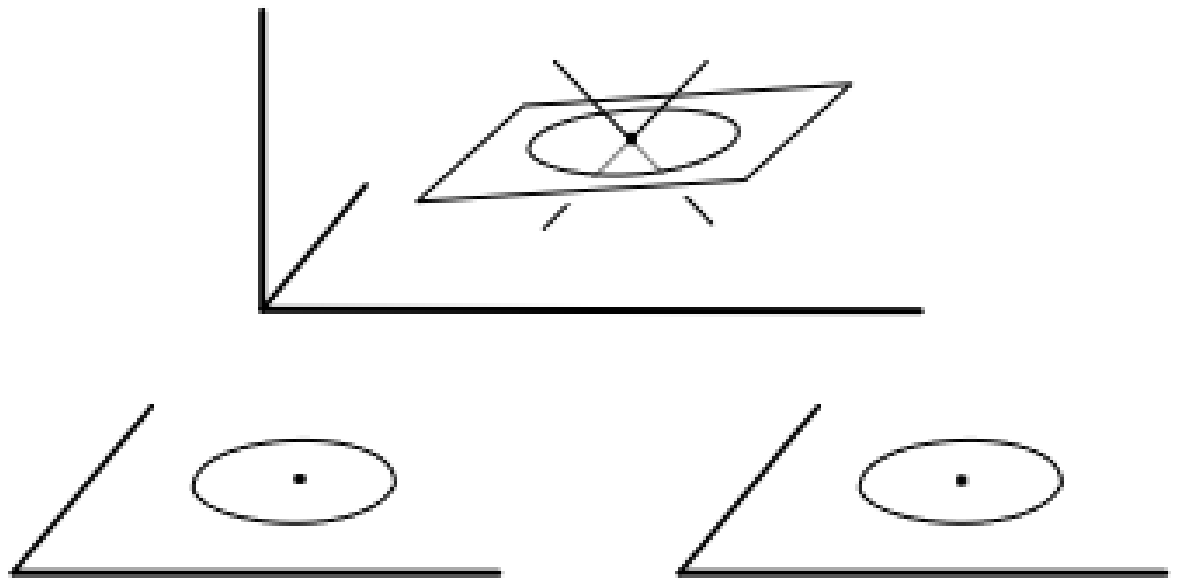} }
 },
 (2.2,7)*\txt{$K_0$},
 (1,2.3)*\txt{$H'_0$},
 (9.2,2.3)*\txt{$H_0$},
 (7.4,7.1)*\txt{$D_m\subset F_m$},
 (5.2,8.2)*\txt{$m$},
 (4.1,2.6)*\txt{$U'_m$},
 (8.4,2.6)*\txt{$U_m$},
 (2.8,2.3)*\txt{$x'$},
 (7.6,2.25)*\txt{$x$},
 (3.1,3.8)*\txt{$\theta_0$},
 (7.3,3.7)*\txt{$\eta_0$},
 \ar (3.99,4.4) ; (3.36,3),
 \ar (6.41,4.4) ; (7.04,3)
\endxy
}
The disc $D_m$ is transversal to the $\eta_0$- and $\theta_0$-fibers at $m$, thus we may assume (shrinking if necessary) that the restrictions of $\eta_0$ and $\theta_0$ to $D_m$ are diffeomorphisms onto open sets $U_m\subset H_0$ and $U'_m\subset H'_0$. We may also assume that there exists a chart around $m$, containing $D_m$ in its domain, in which the submersion $\eta_0$ is locally the projection of a product $H_0\times F\to H_0$.\\
Set $Z_0$ the disjoint union of such discs $D_m$ for all $m\in K_0$. We claim that the canonical map $j:Z_0\to K_0$ induces an essential equivalence. We first check that $t\pi_K:K_s\times_j Z_0\to K_0$ is a surjective submersion by constructing local sections of this map (see figure on the next page).\\[.3em]
Take any $n\in K_0$, and any $(k,m)\in K\times Z_0$ in the $t\pi_K$-fiber over $n$, that is: \f{k}{m}{n} in $K$. Denote $\Db$ the reference disk for $m$ in $Z_0$. There exists a local bisection $\beta$ extending $k^{-1}$ in a neighborhood $V$ of $n$ [\ref{lbeta}], we shrink it if necessary to have $\eta t\beta(V)\subset U_{m^*}$.
We then crush $W=t\beta(V)$ into $\Db$ following the arrows given by
$$\xymatrix @W=56mm @R=-2pt @C=7mm {
*+=<10mm,0mm>{W}
 \ar [r]
 \ar @{} [d] _{\ds p\ :\quad}
& K_W^{\Db}=s^{-1}(W)\cap t^{-1}(\Db)\\
*+=<10mm,0mm>{m'}
 \ar@{|->} [r]
& (t,\eta,s)^{-1}\big(\,\sigma\eta (m')\,,\,1_{\eta m'}\,,\,m'\,\big)
}$$
where $\sigma$ is the inverse of $\eta_0:\Db\to U_{m^*}$. Following $p^{-1}$ and then $\beta^{-1}$, we get a section of $t\pi_K$, which may be written precisely $(\beta^{-1}\cdot (pt\beta)^{-1}, j^{-1}tpt\beta)$ (where $j$ is restricted to $\Db$ for the inverse).\\
\centerline{
\xy
\xyimport(10,10){
 \resizebox{!}{55mm}{ \includegraphics{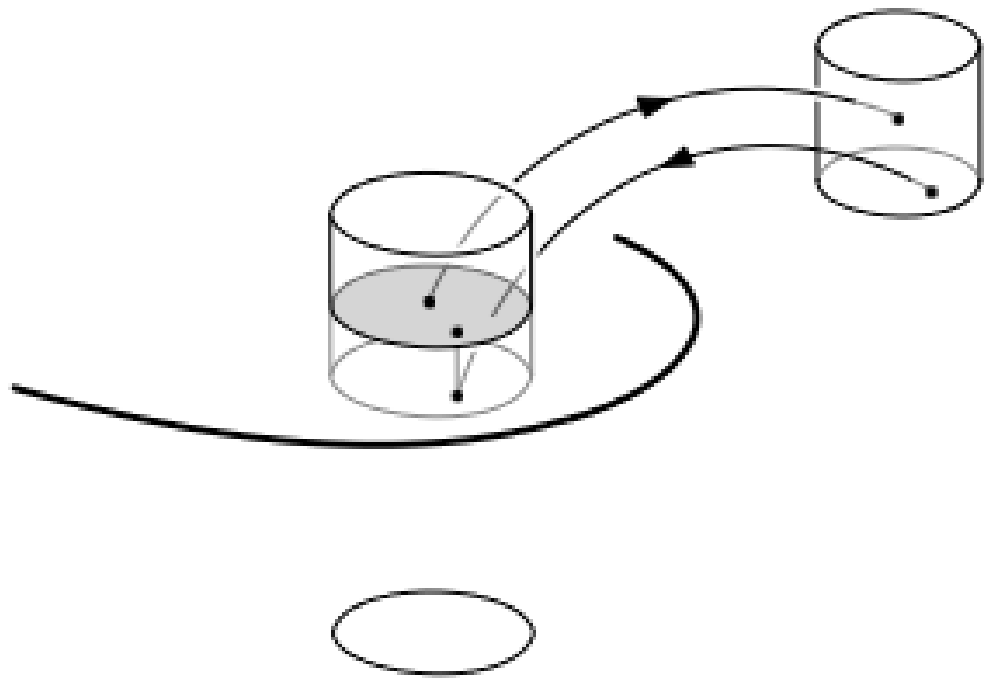} }
 },
 (3.25,7.2)*\txt{$W$},
 (8.4,9.3)*\txt{$V$},
 (5.7,8.47)*\txt{$k$},
 (6.3,6.7)*\txt{$\beta n'$},
 (7.8,7.8)*\txt{$n$},
 (8.1,6.5)*\txt{$n'$},
 (3.96,5.47)*\txt{$m$},
 (4.25,4.5)*\txt{\gc $p$},
 (1.8,5)*\txt{$\Db$},
 (3.9,2.7)*\txt{$\sigma$},
 (4.8,2.7)*\txt{$\eta_0$},
 \ar @^{->} (4.4,3.4) ; (4.4,2.1),
 \ar @^{->} (4.3,2.1) ; (4.3,3.4),
\endxy
}
It follows that $j$ induces a groupoid pullback [\ref{pbgd}], with an essential equivalence for canonical map $J:Z\to K$. Let us denote $\phi=\eta J$ and $\psi=\theta J$. We claim that $\phi,\psi:Z\to H,H'$ fit the problem.\\
$$\xymatrix {
& Z
 \ar [dl] _{\psi}
 \ar [dr] ^{\phi}
 \ar @[lagray] [d] ^(.52){\gc J} \\
H'
&{\gc K}
 \ar @[lagray] [l] ^{\gc\theta}
 \ar @[lagray] [r] _{\gc\eta}
& H
}$$
By definition of $Z_0$, we know that it has dimension $d$, and that $\phi_0$ and $\psi_0$ are surjective and \'etale. Thus it only remains to check whether the canonical functors $Z\to\phi_0^*H$ and $Z\to\psi_0^*H'$ are diffeomorphisms. This can be achieved by writing $Z=Z_0\,{}_j\times_t K_s\times_j Z_0=Z_0\,{}_j\times_{\id}K_0\,{}_{\eta_0}\times_t H_s \times_{\eta_0} K_0\,{}_{\id}\times_j Z_0=Z_0\,{}_{(\eta_0 j)}\times_t H_s \times_{(\eta_0 j)}Z_0$, and same with $H'$. \cqfd

\main\label{m2hproof}{\em Proof of proposition \ref{mor2haef}.} Considering proposition \ref{mordim}, it suffices to prove the statement when we have a Morita morphism $\phi:H\to H'$. In this case, $\phi_0:H_0\to H'_0$ is a surjective submersion between manifolds of the same dimension, thus it is a surjective local diffeomorphism. Cover $H_0$ with open sets $U_i$ such that the restrictions $\phi_i=\phi_{0|U_i}:U_i\to V_i$ are diffeomorphisms, and set $\Phi=\{\phi_i\}$. We claim that $\Phi$ gives a Haefliger equivalence from $\tloc H$ to $\tloc H'$.\\
Due to the stability of pseudogroups under gluing, it is enough to prove the following: for any bisection $\beta:U\to H_{U_i}^{U_j}$ with domain contained in a single $U_i$ and $t\beta(U)$ contained in a single $U_j$ (resp $\beta':V\to(H')^{V_j}_{V_i}$), the composition $\phi_j\circ t\beta\circ\phi_i^{-1}$ is an element of $\tloc H'$ (resp $\phi_j^{-1}\circ t\beta'\circ\phi_i\in\tloc H$). But all such bisections $\beta$ and $\beta'$ are in one-to-one correspondance through:\\
\centerline{
\xy
\xyimport(10,10){
 \resizebox{!}{40mm}{ \includegraphics{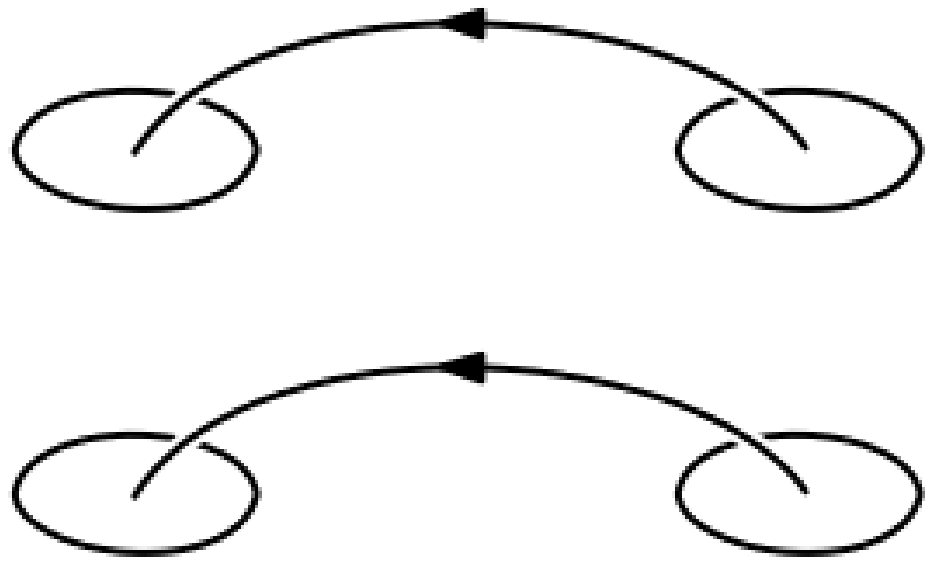} }
 },
 (5.2,8.6)*\txt{$\beta$},
 (5.3,2.6)*\txt{$\beta'$},
 (1.35,6.5)*\txt{$H$},
 (1.35,2.3)*\txt{$H'$},
 (2.2,4.2)*\txt{$\phi_j$},
 (8.3,4.2)*\txt{$\phi_i$},
 (5.7,5.6)*\txt{$\phi$},
 \ar (2.83,5.2) ; (2.83,3.3),
 \ar (7.58,5.2) ; (7.58,3.3),
 \ar (5.2,7) ; (5.2,4.3)
\endxy
}
$$\xymatrix @W=65mm @R=-2pt @C=7mm {
*+=<10mm,0mm>{\beta}
 \ar @{|->} [r]
& \beta'=\phi\beta\phi_i^{-1}\\
*+=<10mm,0mm>{\beta'}
 \ar @{|->} [r]
&\beta=(t,\phi,s)^{-1}(\phi_j^{-1}t\beta'\phi_i,\beta'\phi_i,\id)
}$$
Thus every $\phi_j t\beta\phi_i^{-1}=t\phi\beta\phi_i^{-1}$ is some $t\beta'\in\tloc H'$, and every $\phi_j^{-1}t\beta'\phi_i$ is some $t\beta$, with $\beta$ given by the correspondance. \cqfd
\vspace*{.2em}
\begin{defenv}{[Groupoid dimension]}
 The {\df groupoid dimension} of a Lie groupoid $H$ is the relative integer $\gdim H=\dim H_1-2\dim H_0$.
\end{defenv}
It is immediate from the definition of pullback groupoids [\ref{pbgd}] that the groupoid dimension is preserved under pullbacks, and therefore under Morita-equivalence [\ref{moree}].\\[.3em]
{\em Proof of theorem \ref{thmpggd}.} According to propositions \ref{haef2mor} and \ref{mor2haef}, the only point missing to get the theorem is to check whether the groupoids $H=[\pgH]$ and $H'=[\pgH']$ have bases of the same dimension whenever they are Morita-equivalent. In this case, we know that $\gdim H=\gdim H'$. But $H$ and $H'$ are \'etale, therefore $\gdim H=\dim H_0-2\dim H_0=-\dim H_0$ and $\gdim H'=-\dim H'_0$, and thus $\dim H_0=\dim H'_0$. \cqfd

\section{Compact generation}

\main\label{obsep} In this section we will define a property for groupoids involving compactness. The problem is that we have to deal with fiber products over base spaces, which implies heavy manipulations on compact subspaces that may go wild in the non-Hausdorff case. Therefore we will restrict ourselves to a (not so) special kind of groupoids:
\pagebreak
\begin{defenv}{}
 A Lie groupoid is called {\bf object-separated} if its base is a Hausdorff space.
\end{defenv}
Groupoids in classical foliation theory are naturally object-separated, as their bases are (smooth) Hausdorff manifolds. {\em All groupoids will be hereafter assumed object-separated.} This definition is naturally compatible with pullbacks to Hausdorff bases, so that we can go on using the constructions defined in the first section, however we need to refine Morita equivalence [\ref{moree},\ref{mormor}] a bit to fit our new type of groupoids:
\begin{prop}{}\label{morsep}
 Let $G$ and $H$ be two object-separated, Morita-equivalent groupoids. Then there exists an object-separated groupoid $K$, and two Morita morphisms $\phi,\psi:K\to G,H$.
\end{prop}
{\em Proof: `covering trick'}. Let $K'$ be a groupoid yielding a Morita equivalence between $G$ and $H$, with two Morita morphisms $\phi',\psi':K'\to G,H$ (theorem \ref{mormor}). Cover $K'_0$ with open sets $U_i$ diffeomorphic to $\R^n$, and let $K_0$ be the disjoint union of the $U_i$'s (which is a Hausdorff manifold). The canonical map $j:K_0\to K'_0$ is a surjective submersion, and therefore induces a groupoid pullback $K:=j^*K'$, with a Morita morphism $J$ as canonical functor (see \ref{pbgd}). Set $\phi:=\phi'\circ J$, $\psi:=\psi'\circ J$, these are Morita morphisms as compositions of Morita morphisms (straightforward exercise). \cqfd

\main Before defining compact generation, let us introduce some vocabulary. Any Lie groupoid $G$ `acts'\footnote{Can be formalized \cite{Haef84}.} naturally on itself by left multiplication. For any $x\in G_0$, we shall call $G\cdot x=t(s^{-1}(x))$ the $0$-orbit of $x$ in $G$~; it may be seen as the set of all points of $G_0$ which are linked to $x$ by an arrow of $G$. In a similar way, we define the {\df $\mathbf 1$-orbit} of $g\in G$ to be the set of all arrows of $G$ which are linked to $g$ by an arrow of $G$, that is $\Orb(g)=t^{-1}(t(s^{-1}(sg)))=s^{-1}st^{-1}tg$.\\
If $S$ is any subset of $G$, we define the {\df base} of $S$ to be $S_0:=\langle S\rangle_0=s(S)\cup t(S)$. We shall say that $S$ is {\df exhaustive in $\mathbf G$} if it meets every $1$-orbit, or equivalently if its base meets every $0$-orbit. Finally, we recall that relative compactness means inclusion into a compact subset [\ref{relk}].
\pagebreak
\begin{defenv}{}\label{gk}
 An object-separated Lie groupoid $G$ is said to be {\df compactly generated} if it contains a relatively compact exhaustive open subset $\gU$, which generates a full subgroupoid.
\end{defenv}
By propostion \ref{openprod}, we know that an open subset of $G$ generates an open subgroupoid, for the generated set $\langle\gU\rangle$ is just the union of all powers of $\gU\cup\gU^{-1}$. In particular, $\langle\gU\rangle$ is a Lie subgroupoid of $G$, and it may be seen that this definition is equivalent to ``there exists a relatively compact open subset $\gU$ such that $\langle\gU\rangle\subset G$ is an essential equivalence''.\\
For convenience, we shall say that a subset $S\subset G$ which generates a full subgroupoid has the full generation property, or simply has full generation.

\main We have seen that pseudogroups and germ groupoids could be naturally identified (theorem \ref{thmpggd}). There is already a notion of compact generation for pseudogroups, therefore we begin by investigating the relation between definition \ref{gk} for germ groupoids and the classical definition for pseudogroups:
\begin{defenv}{}
 A pseudogroup $(\pgH,T)$ is said to be {\bf compactly generated} if the following holds:
 \vspace*{-2em}
 \begin{center}\begin{tabular}{c p{.95\textwidth}}
  $\bullet$&$T$ admits an open subset $U$, relatively compact, and exhaustive (meeting every $\pgH$-orbit).\\
  $\bullet$& There exist finitely many $h_i\in\pgH$ and open sets $U_i\subset U$, each relatively compact in the domain of the associated $h_i$, such that the induced pseudogroup $\pgH_{|U}$ (elements of $\pgH$ with domain and image in $U$) is generated by the $(h_i)_{|U_i}$.
 \end{tabular}\end{center}
\end{defenv}
\begin{prop}{}\label{cgpggd}
 Let $(\pgH,T)$ be a pseudogroup and $H$ its germ groupoid. Then $H$ is compactly generated if and only if $\pgH$ is.
\end{prop}
\proof Assume that $\pgH$ is compactly generated, and consider the set of germs
$$\gU=\cups_i\left\{[h_i]_x\;;\;x\in U_i\right\}$$
This is an open subset of $[\pgH]$ [\ref{pg2gd}], relatively compact as finite union of relatively compact sets (check it in the charts provided by the domains of the $h_i$'s), and exhaustive because $U$ is. It also generates a full subgroupoid because every arrow of $[\pgH]$ between points of $\bgU=U$ is a germ of some element in $\pgH_{|U}$, so is a product of germs of the $(h_i)_{|V_i}$ at suitably chosen points. Therefore $[\pgH]$ is compactly generated.\\
Conversly, assume that $H$ is compactly generated, and let $\gU$ be a symmetric exhaustive relatively compact open subset of $[\pgH]$ with full generation (we may suppose $\gU$ symmetric, for $\gU\cup\gU^{-1}$ has the same properties as $\gU$ regarding compact generation). Decompose $\gU$ into a finite union of open subsets $\gU_i$ such that each $\gU_i$ is included in a compact set $\gK_i$, itself included in an open set $\gV_i$ where $s$ and $t$ are both diffeomorphisms onto open sets of $T$ (the union is finite because $\gU$ is included in a compact). Write $s_i=s_{|\gV_i}$, and set $h_i=t\circ s_i^{-1}\in\tloc[\pgH]=\pgH$, $U_i=s(\gU_i)$ and $U=\gU_0=\cups_i U_i$. We claim that this data satisfies the definition of compact generation for $\pgH$. It is not hard to see that $U$ is an exhaustive relatively compact open subset of $T$: exhaustiveness is inherited from that of $\gU$ ; openness and relative compactness are consequences of $s$ and $t$ being \'etale (and $\gU$ relatively compact). Now take some $h\in\pgH_{|U}$, and any $x$ in the domain of $h$. The germ $[h]_x$ has source and target in $\gU_0$, hence by full generation it may be written as a product of elements of $\gU$:
$$[h]_x = u_l\cdot\ldots\cdot u_1$$
Each $u_k$ is in some $\gU_{\alpha(k)}$, and the definition of $[\pgH]$ implies $u_k=[h_{\alpha(k)}]_{su_k}$. Thus
$$[h]_x = [h_{\alpha(l)}]_{su_l}\cdot\ldots\cdot [h_{\alpha(1)}]_{su_1} = [h_{\alpha(l)}\ldots h_{\alpha(1)}]_x$$
and $h$ is locally at $x$ a product of the $h_{i|U_i}$'s. As it is true for every $x$ in its domain, by gluing, $h$ is in the pseudogroup generated by the $h_{i|U_i}$'s. \cqfd

\main {\sc Example.}\label{example} Given a foliation $\fol$ on a manifold $M$, a classical groupoid associated to $\fol$ is its holonomy groupoid $\Hol(\fol)$, the set of all germs of holonomy elements\footnote{Recall that a holonomy element is a diffeomorphism associated to a tangent path, which roughly follows the transverse coordinates of the neighboring leaves along the path.} of the foliation up to a choice of local transversals. This groupoid is the modern evolution of the holonomy pseudogroup $\pgH_T$ associated to a complete transversal $T$ \cite{gdbfol}. A classical result asserts that the germ groupoid of this pseudogroup identifies to the pullback of $\Hol(\fol)$ along the inclusion $T\to M$, the resulting morphism being an essential equivalence [\ref{ee}]. The same pullback construction may be achieved for the other classical groupoid associated to $\fol$, the groupoid of tangent paths up to tangent homotopy - called the monodromy (or homotopy) groupoid $\Mon(\fol)$ of $\fol$, producing another \'etale groupoid. In particular, these two classical groupoids are Morita-equivalent to \'etale ones. These remarks have lead to the definition of `foliation' groupoids \cite{MoerdijkFolGd}:
\begin{defenv}{}
 A foliation groupoid is any groupoid which is Morita-equivalent to an \'etale one.
\end{defenv}
The $0$-orbits of a foliation groupoid naturally define a foliation of its base \cite{MoerdijkFolGd}. Given that compact generation was originaly designed to characterize holonomy pseudogroups of compact foliated manifolds, we ask the question: is a foliation groupoid with compact base always compactly generated?\\
\begin{thmstyle}{\sc Theorem \notmain}\label{cg4cb}
 Every $s$-connected foliation groupoid with compact base is compactly generated.
\end{thmstyle}
Recall that an $s$-connected groupoid is a groupoid with connected $s$-fibers (or $t$-fibers). The theorem fails to be true if we don't require the groupoid to be $s$-connected: choose any non-finitely generated group $A$ and any compact manifold $M$. Give $A$ the discrete topology and consider the trivial groupoid on $M$ with group $A$, $M\times A\times M$ (with product $(z,b,y)(y,a,x)=(z,ba,x)$). It is a foliation groupoid: take any $m\in M$, and check that the inclusion $\{m\}\times A\times\{m\}\subset M\times A\times M$ is an essential equivalence. It has a compact base, but cannot be compactly generated. Indeed, if it was, there would exist a relatively compact subset $\gU\subset M\times A\times M$ with full generation. The group $A$ is not finitely generated, so $\gU$ would have to cross an infinite number of $M\times\{a\}\times M$ to have full generation. But it is impossible, because the sets $M\times\{a\}\times M$ are open and pairwise disjoint, so that $\gU$, which is relatively compact, can only meet a finite number of them.\\[.5em]
We denote $G$ a foliation groupoid, and $\fol$ the foliation it induces on its base, which we assume to be compact. The next proof uses both the local structure of foliation groupoids and the construction of the natural factorisation morphsim $h_G:\Mon(\fol)\to G$ given in \cite{MoerdijkFolGd}, which we recall here.
\begin{thmstyle}{\\\cite{MoerdijkFolGd} Lemma 3}
 A trivializing submersion $\pi:U\to T$ of $\fol$ is a submersion on an open set of $G_0$ with contractible fibers, which are precisely the leaves of $\fol_{|U}$.\\
 Denote $G(U)$ the $s$-connected component of $G_U^U$. Then the map
 $(t,s):G(U)\to U\times_T U$
 is a natural isomorphism of Lie groupoids, where $U\times_T U$ is given the pair product $(z,y)(y,x)=(z,x)$.
\end{thmstyle}

\pagebreak\nin
\begin{thmstyle}{\cite{MoerdijkFolGd} Proposition 1}
 With the preceding notations, there is a natural factorisation of the canonical map $hol:\Mon(\fol)\to\Hol(\fol)$ into two surjective (functorial) local diffeomorphisms
 \vspace*{-0.3em}
 $$\xymatrix@1 @M +1mm { \Mon(\fol)\ar[r]^(.6){h_G} & G \ar[r]^(.4){hol_G} & \Hol(\fol) } $$
\end{thmstyle}
{\em Proof of thm \ref{cg4cb}.} Let $(V_i)_i$ be a finite open cover of $G_0$ by domains of trivializing submersions, and $(U_i)_i$ a shrinking of $(V_i)_i$. Set $\gU$ to be the union of the restrictions $G(U_i)=G(V_i)_{U_i}^{U_i}$ of the local groupoids $G(V_i)$. The set $\gU$ is open in $G$, and clearly exhaustive because $\bgU=G_0$. Using the isomorphism $(t,s)_i:G(V_i)\to V_i\times_{T_i}V_i$, we can include each $G(U_i)$ in a compact set (namely $(t,s)_i^{-1}\big(\overline{U_i}\times_{T_i}\overline{U_i}\big)$), which implies that $\gU$ is relatively compact. Finally, as $(U_i)_i$ covers $G_0$, every tangent path $\alpha$ of $G_0$ may be decomposed as a product of tangent paths $\alpha_k$ contained each in a single $U_{i(k)}$. Due to the contractibility of the $\pi_{i(k)}$-fibers, each path $\alpha_k$ is entirely defined by its ends. Those ends in turn give an element $g_k=(t,s)_{i(k)}^{-1}(\alpha_k)\in G(U_{i(k)})$. Then the product of the $g_k$'s lies in the subgroupoid generated by $\gU$, and is precisely $h_G(\alpha)$ by construction of $h_G$. It follows that $\gU$ generates the entire image of $h_G$. But $G$ is $s$-connected, so that $h_G$ is surjective onto $G$. Thus $\gU$ generates $G$, and has full generation. \cqfd
\begin{cor}{}\label{cg4folg}
 If $\fol$ is a foliation on a compact manifold, its monodromy and holonomy groupoids are compactly generated.
\end{cor}
\proof It suffices to remark that the $s$-fibers of a monodromy (resp. holonomy) groupoid are the universal (resp. holonomy) coverings of the leaves, and thus are connected. \cqfd

\begin{mainthm}{}\label{invgk}
 If $G$ and $H$ are Morita-equivalent object-separated Lie groupoids, and if $G$ is compactly generated, then so is $H$.
\end{mainthm}
Considering proposition \ref{morsep}, it suffices to prove the following
\begin{lemma}{}
 If $\phi:G\to H$ is a Morita morphism between object-separated groupoids, and if $G$ or $H$ is compactly generated, then so is the other one.
\end{lemma}
The easy case is when $G$ is compactly generated: let $\gU$ be a relatively compact exhaustive open subset of $G$ with full generation, and let $\gV:=\phi(\gU)$. Then $\gV$ is immediately relatively compact and open, as $\phi$ is a submersion (check in the complete pullback diagram \ref{pbgd}, with $f=\phi_0$ a submersion). The unique lifting property [\ref{ee}] shows that the $1$-orbits of $G$ and $H$ are in one-to-one correspondance via $\phi$ ($\Orb(g)\mapsto\Orb(\phi g)$), so that $\gV$ meets every $H_1$-orbit which is the image of a $G_1$-orbit met by $\gU$. But $\gU$ meets every of them, thus $\gV$ is exhaustive. Finally, as $\phi$ is a functor, the subgroupoid generated by $\gV$, the image of $\gU$, is the image of the subgroupoid generated by $\gU$, which is full. Thus $\langle\gV\rangle$ is also full, as the unique lifting property ensures it can't miss any arrow.\\[.6em]
The other case is a bit tricky, because we have to climb up $\phi$ preserving both openness and relative compactness. Assume we have a $\gV\subset H$ giving compact generation for $H$, chosen symmetric.\\
The map $\phi_0$ is a submersion, so we can locally write it as the projection of a product onto an open set $W\times F\to W$. As $V:=\bgV$ is relatively compact, we can cover it with a finite number of open sets $V_i\subset V$, each included in a compact $K_i$, itself contained in the image of a trivialization $W_i\times F_i\to W_i$ of $\phi_0$. For each $i$, we also choose a small open set $D_i$, inside a compact subset $F'_i\subset F_i$, and set $U_i:=V_i\times D_i$, $C_i:=K_i\times F'_i$. Let us write $U:=\cups_i U_i$, and consider the set of all arrows between points of $U$ that are sent in $\gV$ by $\phi$:
$$ \gU := (t,\phi,s)^{-1}\big(\,U\times\gV\times U\,\big) $$
This is an open subset of $G$, which is included in the compact set
$$ \gC := (t,\phi,s)^{-1}\big(\,C\times\gK\times C\,\big) $$
where $C$ is the union of the $C_i$'s, and $\gK$ some fixed compact set containing $\gV$ ($\gV$ relatively compact). The set $\gC$ is indeed compact, for it is the direct image under the diffeomorphism $(t,\phi,s)^{-1}:\phi_0^*H\to G$ of the intersection of the compact cartesian product $C\times\gK\times C$, and the {\em closed} submanifold $\phi_0^*H\subset G_0\times H\times G_0$ (here we use that $\phi_0^*H=(G_{0\,\phi_0})H(_{\phi_0}G_0)$ is a fiber product over twice the diagonal of $G_0$, which is {\em closed} [\ref{obsep}]). Thus $\gU$ is relatively compact.\\
It is easy to see that $\gU$ is exhaustive: by unique lifting every $G_1$-orbit is the preimage by $\phi$ of an $H_1$-orbit, which necessarily crosses $\gV$ (exhaustiveness), and thus its preimage crosses $\gU$ ($\phi(\gU)=\gV$ by construction).\\
It only remains to check the full generation property. Take any $g\in G_U^U$. We have $\phi g\in H_V^V=\langle\gV\rangle$, thus we can write $\phi g$ as a finite product of elements of $\gV$
$$ \phi g = v_l\cdot\ldots\cdot v_1 $$
Choose indices $\alpha(k)$ such that $x_0:=sv_1\in V_{\alpha(0)}$, and $x_k:=tv_k\in V_{\alpha(k)}$ for all $k>0$. We then choose points $m_k\in U_{\alpha(k)}$ over the $x_k$'s by $\phi_0$, with two special choices $m_0=sg$ and $m_l=tg$, and define
$$ g_k:=(t,\phi,s)^{-1}\big(\,m_k\,,\,v_k\,,\,m_{k-1}\,\big)
\quad\text{for }0<k<l$$
\centerline{
\xy
\xyimport(10,10){
 \resizebox{!}{50mm}{ \includegraphics{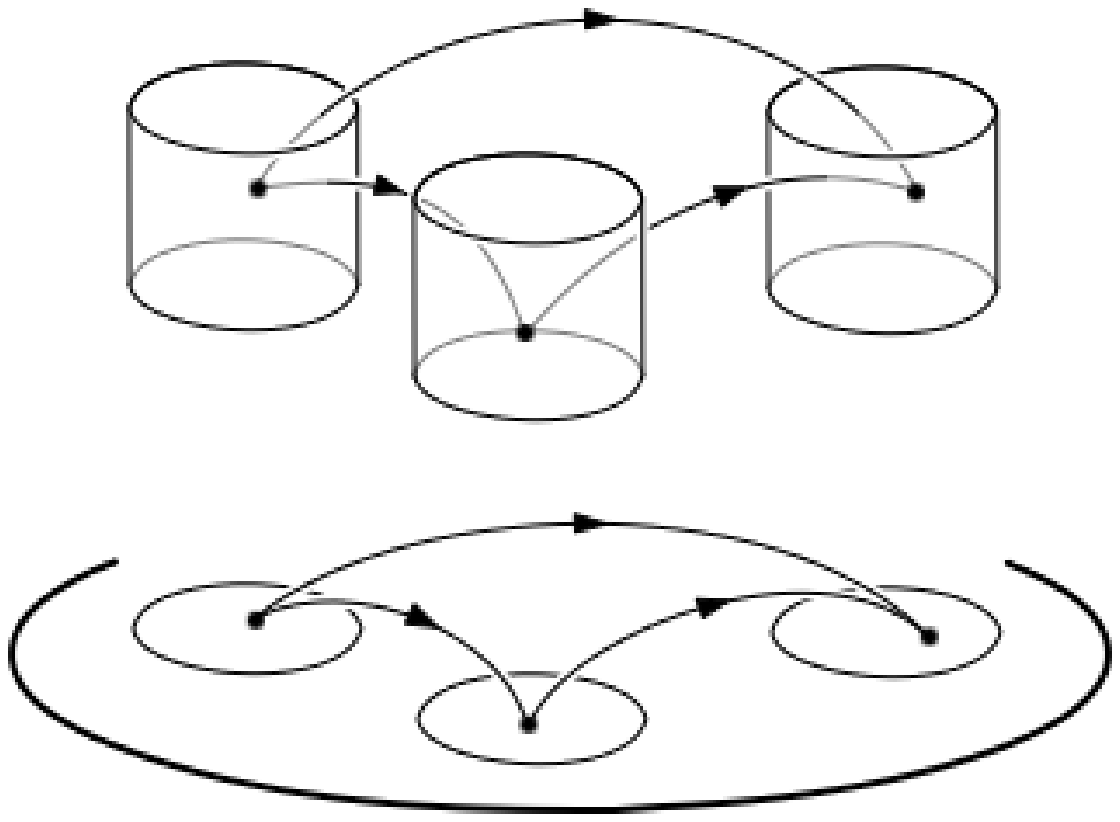} }
 },
 (5.8,9.4)*\txt{$g$},
 (1.75,7.2)*\txt{$sg$},
 (8.3,7.2)*\txt{$tg$},
 (4.9,5.2)*\txt{$m_1$},
 (5.8,4.2)*\txt{$\phi g$},
 (3.92,2.2)*\txt{$v_1$},
 (6.2,2.2)*\txt{$v_2$},
 (.4,7)*\txt{$G_0$},
 (.4,2.3)*\txt{$V$},
 (-.1,4.6)*\txt{$\phi_0$},
 \ar (.4,6.2) ; (.4,3)
\endxy
}
By definition of $\gU$, $g_k\in\gU$ for all $j$, so that the (well-defined) product $g_l\ldots g_1$ is in $\langle\gU\rangle$. But
$$ \phi(g_l\ldots g_1)\eq\phi g_l\cdot\ldots\cdot\phi g_1\eq v_l\cdot\ldots\cdot v_1\eq\phi g $$
and unicity of the pullback between $m_0=sg$ and $m_l=tg$ implies $g=g_l\ldots g_1\in\langle\gU\rangle$. Thus $\gU$ has full generation, and $G$ is compactly generated. \cqfd
\begin{cor}{[Haefliger's lemma, groupoid version]}
 If $\fol$ is a foliation on a compact manifold $M$, all transverse holonomy groupoids of $\fol$ are compactly generated.
\end{cor}
We have seen (theorem \ref{thmpggd} and proposition \ref{cgpggd}) that this assertion is simply a translation of Haefliger's original result in groupoid theory.\\
\proof We know by corollary \ref{cg4folg} that $\Hol(\fol)$ is compactly generated. Given a complete transversal $T$ for $\fol$, it is a classical result that the inclusion $T\subset M$ induces an essential equivalence $\Hol_T(\fol)\to\Hol(\fol)$ [\ref{ee}], where $\Hol_T(\fol)$ is the germ groupoid of the holonomy pseudogroup $\pgH_T$ associated to $T$. In particular these two groupoids are Morita-equivalent [\ref{moree}], and thus $\Hol_T(\fol)$ is compactly generated by the invariance theorem \ref{invgk}. \cqfd

\main {\sc Remark.} As in the case of pseudogroups, it is important to require $\gU$ to be open in definition \ref{gk}. We construct here an example of a non-compactly generated groupoid, which admits a relatively compact exhaustive non-open subset $\gU$ with full generation. This example is inspired from an exercise in \cite{Haef00}.\\
Consider the Klein bottle, seen as a cylinder $S^1\times[-1,1]$ ($S^1\subset\C$) with its ends identified $(z,1)=(\overline{z},-1)$, and foliated by the directrices $\{z\}\times[-1,1]$.\\
\centerline{
\xy
\xyimport(10,10){
 \resizebox{!}{50mm}{ \includegraphics{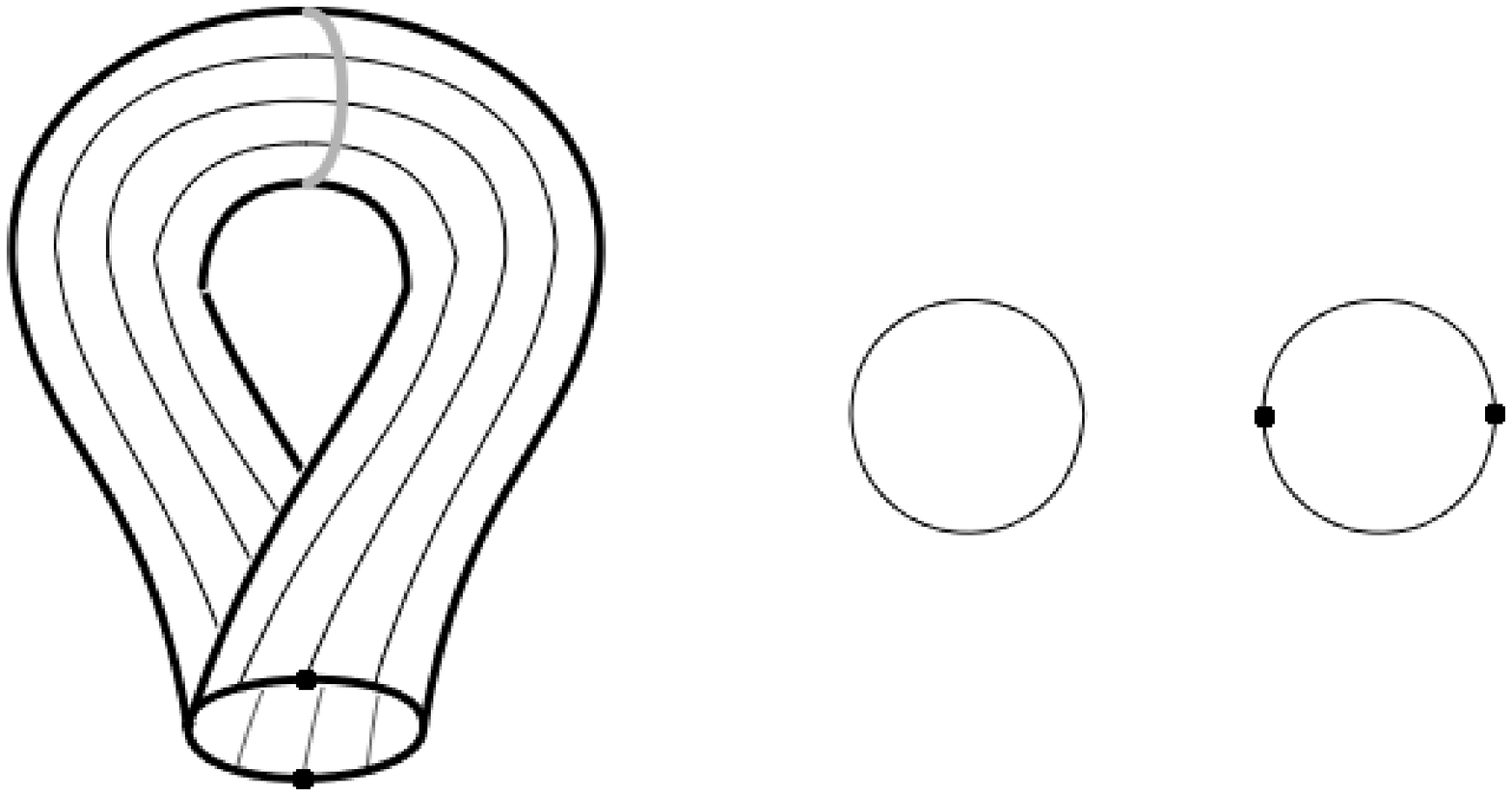} }
 },
 (3.1,9.1)*\txt{\gc $T$},
 (.6,4.8)*\txt{$M$},
 (6.13,4.3)*\txt{$[\id]$},
 (8.37,4.35)*\txt{$[\conj]$},
 (7.25,4.3)*\txt{$\coprod$},
 (7.15,6.5)*\txt{$\Hol_T(\fol)=$}
\endxy
}
Take $T=S^1\times\{0\}$ a circular transversal, then the transverse holonomy groupoid $H=\Hol_T(\fol)$ associated to $T$ may be seen as two copies of $S^1$, one for trivial holonomy germs and one for the germs of the holonomy diffeomorphism obtained by turning once along the directrices (which induces the complex conjugation on $S^1$).\\
Now make two holes in the Klein bottle at $(\pm 1,1)$, so that we cannot turn around the associated leaves. In this case, $H$ loses the points $[\conj]_{\pm 1}$, and is no more compactly generated. Indeed, if we had some $\gU$ satisfying definition \ref{gk}, $\bgU$ would contain some small neighborhood $V$ of $1$, which we may assume stable under conjugation. The full generation property would then force $\{[\conj]_z ; z\in V\setminus\{1\}\}\subset\gU$, and thus $\gU$ could not be relatively compact in $H$, absurd.
However, if we set
$$\gU=\{[\id]_z ; \Im m(z)\geqslant 0 \}$$
then $\gU$ is relatively compact, exhaustive, and generates a full subgroupoid of
$H$. It statisfies all conditions for compact generation, {\em but} being open.

\nocite{PCartier,Haef70,Haef84,Haef00,MackenzieI,MackenzieII,MoerdijkFolGd,MrcunKT,MrcunPhD,PingXu}
\bibliographystyle{alpha}
\bibliography{biblio_cg4gd}

\begin{thebibliography}{Mac05}

\bibitem[BX06]{PingXu}
Kai Behrend and Ping Xu.
\newblock Differentiable stacks and gerbes, 2006.
\newblock
  \href{http://arxiv.org/abs/math/0605694}{\url{http://arxiv.org/abs/math/0605%
694}}.

\bibitem[Car08]{PCartier}
Pierre Cartier.
\newblock Groupo\"ides de lie et leurs alg\'ebro\"ides.
\newblock {\em S\'eminaire Bourbaki}, 2007-2008:987.1--987.30, 2008.

\bibitem[CM00]{MoerdijkFolGd}
Marius Crainic and Ieke Moerdijk.
\newblock Foliation groupoids and their cyclic homology, 2000.
\newblock
  \href{http://arxiv.org/abs/math/0003119}{\url{http://arxiv.org/abs/math/0003%
119}}.

\bibitem[God91]{gdbfol}
Claude Godbillon.
\newblock {\em Feuilletages, \'Etudes g\'eom\'etriques}.
\newblock Birkhäuser, 1991.

\bibitem[Hae70]{Haef70}
Andr\'e Haefliger.
\newblock Homotopy and integrability.
\newblock In {\em Manifolds}, volume 197 of {\em Lecture Notes in Mathematics},
  pages 133--163. Springer, 1970.

\bibitem[Hae84]{Haef84}
Andr\'e Haefliger.
\newblock Groupo\"ides d'holonomie et classifiants.
\newblock In {\em Structure transverse des feuilletages}, volume 116 of {\em
  Ast\'erisque}, pages 70--97. Soci\'et\'e Math\'ematique de France, 1984.

\bibitem[Hae00]{Haef00}
Andr\'e Haefliger.
\newblock Foliations and compactly generated pseudogroups.
\newblock In {\em Foliations : Geometry and Dynamics}, pages 275--295. World
  Scientific, 2000.

\bibitem[Mac87]{MackenzieI}
Kirill~C.H. Mackenzie.
\newblock {\em Lie Groupoids and Lie Algebroids in Differential Geometry}.
\newblock Cambridge University Press, 1987.

\bibitem[Mac05]{MackenzieII}
Kirill~C.H. Mackenzie.
\newblock {\em General Theory of Lie Groupoids and Lie Algebroids}.
\newblock Cambridge University Press, 2005.

\bibitem[Mr{\v c}96]{MrcunPhD}
Janez Mr{\v c}un.
\newblock {\em Stability and invariants of Hilsum-Skandalis maps}.
\newblock PhD thesis, Universiteit Utrecht, 1996.
\newblock
  \href{http://arxiv.org/abs/math/0506484}{\url{http://arxiv.org/abs/math/0506%
484}}.

\bibitem[Mr{\v c}99]{MrcunKT}
Janez Mr{\v c}un.
\newblock Functoriality of the bimodule associated to a hilsum-skandalis map.
\newblock {\em K-Theory}, 18:235--253, 1999.

\end{thebibliography}

\end{document}